\documentclass[3p,sort&compress,times]{elsarticle}

\usepackage{hyperref}

\hypersetup{hidelinks=true}

\usepackage{amsmath}
\usepackage{amssymb}
\usepackage{booktabs}

\usepackage{dcolumn}
	\newcolumntype{d}[1]{D{.}{.}{#1}}
	\newcolumntype{x}[1]{D{x}{\times}{#1}}

\usepackage{graphicx}

\usepackage{mathtools}
	\mathtoolsset{showonlyrefs=true}

\usepackage{siunitx}
\usepackage{subcaption}
\usepackage{xspace}

\providecommand{\aset}[2]{\ensuremath{\mathcal{R}_{#1}^{(#2)}}\xspace}
\providecommand{\bsizes}{\ensuremath{\mathcal{B}}\xspace}
\providecommand{\cset}[2]{\ensuremath{\mathcal{C}_{#1}^{#2}}\xspace}

\providecommand{\elll}{\ensuremath{\ell}}
\providecommand{\kmax}{\ensuremath{K}\xspace}
\providecommand{\kset}[2]{\ensuremath{\mathcal{K}_{#1}^{(#2)}}\xspace}
\providecommand{\lmax}{\ensuremath{L}\xspace}
\providecommand{\mmf}{\ensuremath{\mathrm{MMF}_{\boxplus}}\xspace}
\providecommand{\mmfcsr}{\ensuremath{\mathrm{MMF}_{\mathrm{CSR32}}}\xspace}
\providecommand{\mmflb}{\ensuremath{\mathrm{MMF}_{\mathrm{lb}}}\xspace}
\providecommand{\nnz}{\ensuremath{\mathit{nnz}}\xspace}
\providecommand{\prnnz}{\ensuremath{\mathit{prnnz}}\xspace}
\providecommand{\rnnz}{\ensuremath{\mathit{rnnz}}\xspace}
\providecommand{\schemes}{\ensuremath{\mathcal{S}}\xspace}
\providecommand{\scheme}{\ensuremath{s}\xspace}
\providecommand{\uset}[3]{\ensuremath{\mathcal{U}_{#1,#2}^{#3}}}
\providecommand{\vset}[5]{\ensuremath{\mathcal{V}_{#1,#2,#4}^{#3,(#5)}}}

\DeclareMathOperator{\avg}{avg}

\journal{Computers \& Mathematics with Applications}

\begin{document}

\begin{frontmatter}

\title{On Memory Footprints of Partitioned Sparse Matrices}

\author{Daniel Langr}
\ead{daniel.langr@fit.cvut.cz}
\address{Czech Technical University in Prague\\
Faculty of Information Technology\\
Department of Computer Systems\\
Th\'{a}kurova 9, 160 00, Praha 6, Czech Republic}

\begin{abstract}

Runtime characteristics of sparse matrix computations and related processes may be often improved by reducing memory footprints of involved matrices. Such a reduction can be usually achieved when matrices are processed in a block-wise manner. The presented study analysed memory footprints of 563 representative benchmark sparse matrices with respect to their partitioning into uniformly-sized blocks. Different block sizes and different ways of storing blocks in memory were considered and statistically evaluated. Memory footprints of partitioned matrices were additionally compared with lower bounds and with the CSR storage format. The average measured memory savings against CSR in case of single and double precision were $42.3$ and $28.7$ percents, the corresponding worst-case savings $25.5$ and $17.1$ percents. Moreover, memory footprints of partitioned matrices were in average 5 times closer to their lower bounds than CSR. Based on the obtained results, generic suggestions for efficient partitioning and storage of sparse matrices in a computer memory are provided.


\end{abstract}


\end{frontmatter}


\section{Introduction}
\label{sec:intro}

The way how sparse matrices are stored in a computer memory may have a significant impact on the required memory space, i.e., on matrix memory footprints. Reduction of matrix memory footprints may positively influence related computations and executions of corresponding programs. For example:
\begin{itemize}[---]
\item Lower matrix memory footprints yield faster processing of matrices by I/O subsystems, e.g., when checkpointing-restart resilience methods are applied within high performance computing (HPC) applications~\cite[Chap.~4]{Langr:2014d}\cite{Langr:2013a}. 

\item On modern HPC architectures, the performance of some common sparse matrix operations is highly bounded by memory bandwidth. For instance, during sparse matrix-vector multiplication (SpMV), floating-point units are typically utilized to less than 5 percent of their peak computational capabilities~\cite{Langr:2016a,Williams:2009b}. Lower matrix memory footprints may thus potentially increase the efficiency and performance of sparse matrix computations.

\item Lower matrix memory footprints allow larger matrices to fit in the available amount of memory, which, therefore, allows to solve computational problems to higher extent or with higher accuracy.

\end{itemize}

One way of reducing memory footprints of sparse matrices is their partitioning into blocks (which also promotes spatial locality during computations). Much has been written about block processing of sparse matrices, frequently in the context of memory-bounded character of SpMV; see, e.g., \cite{Belgin:2009,Belgin:2011,Blelloch:1993,Buluc:2009,Buluc:2011,Buono:2016,Byun:2012,Choi:2010,Eberhardt:2016,Im:2001,Im:2004,Kannan:2013,Karakasis:2009a,Langr:2012b,Langr:2016a,Nishtala:2004,Nishtala:2007,Simecek:2012a,Simecek:2015b,Smailbegovic:2005,Stathis:2003,Tvrdik:2006,Williams:2009b}.
In this article, we address the problem of minimizing memory footprints of sparse matrices by their partitioning into uniformly-sized blocks. Its solution raises two essential questions: How to choose a suitable block size? And, how to store resulting nonzero blocks in a computer memory? These questions form a multi-dimensional optimization problem that needs to be solved prior to the partitioning itself. We refer to both these problems---optimization and partitioning---as (\emph{block}) \emph{preprocessing}.

The above introduced optimization problem raises another question: How to specify an optimization space, i.e., a space of tested configurations? Intuitively, the larger the optimization space is, the lower matrix memory footprint may be potentially found, however, at the price of longer preprocessing runtime. To amortize block processing of a sparse matrix, the optimization space thus need to be chosen wisely in a form of a trade-off: we want it to be small enough to ensure its fast exploration but also large enough to contain an optimal or nearly-optimal configuration generally for any sparse matrix.

We present a study that analyses memory footprints of 563 representative sparse matrices from the University of Florida Sparse Matrix Collection (UFSMC)~\cite{Davis:2011} with respect to their partitioning into uniformly sized blocks. These matrices arose from a large variety of applications of multiple problem types and thus have highly diverse structural and numerical properties. Our goal is to minimize memory footprints of matrices and we consider an optimization space that consists of different block sizes and different ways of storing blocks in memory. Based on the obtained results, we finally provide suggestions for both efficient and effective block preprocessing of sparse matrices in general.

\section{Methodology}
\label{sec:meth}

%


In \autoref{sec:intro}, we referred to a \emph{matrix memory footprint} as to an amount of memory space required to store a given matrix in a computer memory. More precisely, we can define it as a number of bits (or bytes) that is needed to store the values of the nonzero elements of a given matrix together with the information about their structure, i.e., their row and column positions.
The ways how sparse matrices are stored in a computer memory are generally called \emph{sparse matrix storage formats}; we call them \emph{formats} only if the context is clear. Matrix memory footprint is thus a function of a given matrix and a used format (memory footprints for the same matrix but distinct formats may differ considerably).

In case of partitioned sparse matrices, their nonzero blocks represent individual submatrices that can be treated separately.
In practice, well-proven formats used for nonzero blocks of sparse matrices are:
%
\begin{itemize}[---]
\item The \emph{coordinate} (COO) format, which stores values of block nonzero elements together with their row and column indices; see, e.g., \cite{Buluc:2009,Langr:2012b,Simecek:2012a}. 

\item The \emph{compressed sparse row} (CSR) format, which stores values and column indices of lexicographically ordered block nonzero elements together with the information about which values / column indices belongs to which block row; see, e.g., \cite{Langr:2012b,Nishtala:2004,Nishtala:2007,Simecek:2012a}. 

\item The \emph{bitmap} format, which stores values of block nonzero elements in some prescribed order and encodes their row and column indices in a bit array; see, e.g., \cite{Buluc:2011,Kannan:2013,Langr:2012b}.

\item The \emph{dense} format, which stores values of both nonzero and zero block elements in a dense array (row and column indices of nonzero elements are thus effectively determined by positions of their values within this array); see, e.g., \cite{Barrett:1994,Im:2001,Im:2004,Langr:2012b}.
\end{itemize}
Considering these formats, we have 6 options how to store nonzero blocks of a sparse matrix in memory:
\begin{enumerate}
\item store {all} the blocks in the COO format,
\item store {all} the blocks in the CSR format,
\item store {all} the blocks in the bitmap format,
\item store {all} the blocks in the dense format,
\item store \emph{all the blocks} in the same format such that the format minimizes the memory footprint of a given matrix (we refer to this option as \emph{min-fixed}),
\item store \emph{each block} generally in a different format such that the format minimizes the contribution of this block to the memory footprint of a given matrix (we refer to this option as \emph{adaptive}).
\end{enumerate}
We call these options \emph{blocking storage schemes}, or shortly \emph{schemes} only. Since the first 4 schemes prescribe a fixed format for all the blocks, we call them \emph{fixed-format schemes}.

For the min-fixed and adaptive schemes, we consider formats for nonzero blocks to be chosen from COO, CSR, bitmap, and dense. In case of the min-fixed scheme, the matrix memory footprint thus contains 2 additional bits for storing the information about the format used for all the nonzero blocks. In case of the adaptive scheme, the matrix memory footprint contains 2 additional bits for each nonzero block to store the information about its format.


To evaluate memory footprints of a given matrix for different schemes and some particular tested block size, we need information about numbers of nonzero elements of all nonzero blocks~\cite{Langr:2012b}. In the end, this information must be obtained for each distinct block size from the optimization space, which represents the most demanding part of the whole optimization process~\cite{Langr:2016b}. The block preprocessing runtime is thus approximately proportional to the number of distinct tested block sizes. Consequently, the lower is their count, the higher are the chances that the partitioning will be profitable at all.


Generally, there is $O(m\times n)$ ways how to set a block size for an $m\times n$ matrix, but for fast block preprocessing, we need to choose only few of them.%
\footnote{In addition to multiplication and Cartesian product, we also use the multiplication sign ``$\times$'' to specify matrix/block sizes. In such cases, $m\times n$ does not denote multiplication, but a matrix/block size of height $m$ and width $n$ (i.e., having $m$ rows ans $n$ columns).}
One possible approach is to consider only block sizes 
\begin{equation}
2^k\times 2^\elll, \quad \text{where}\quad 1\leq k \leq \kmax \quad\text{and}\quad 1\leq \elll\leq \lmax,
\label{eq:kl}
\end{equation}
which reduces the number of tested block sizes to $\kmax\times \lmax$. The rationale behind such a choice is as follows:
\begin{itemize}[---]
\item A core operation for block preprocessing is to find out which nonzero elements belong to which block. 
This operation involves (costly) integer division, however, in case of block sizes $2^k\times 2^\elll$, it may be substituted by much faster logical shift operations. We observed 4 and 7 times faster block preprocessing due to such substitution on an Intel Haswell-based computer system and an Intel Xeon Phi coprocessor, respectively~\cite{Langr:2016b}.

\item COO and CSR formats store local (in-block) row and column indices of block nonzero elements. In case of block sizes $2^k\times 2^\elll$, these indices require exactly $k$ and $\elll$ bits for their memory storage and all these bits are fully employed.

\item Cache lines are typically a power of 2 in size; block sizes $2^k\times 2^\elll$ might thus lead to higher utilization of caches. Consider a multiplication of a sparse matrix by a dense vector (SpMV), which consists of multiplication of matrix nonzero blocks by corresponding vector parts. Assuming a cache line size 64 bytes, block size $8\times 8$, and a double precision computation, exactly one cache line is involved for the input vector as well as for the output vector when multiplying a single nonzero block, provided that vectors are properly aligned in memory. This, among others, allows to develop highly-tuned routines for common block operations (such as their multiplication by dense vector parts), common block sizes, and particular hardware architectures.

\item For the COO, bitmap, and dense formats, nonzero elements of a block of size $2^k\times 2^\ell$ can be stored in so-called \emph{Z-Morton order}~\cite{Morton:1966}, which is based on block's recursive partitioning into sub-blocks of size $2^{k-1}\times 2^{\ell-1}$. This approach may yield higher computational efficiency for both sparse and dense blocks (see, e.g.,~\cite{Buluc:2009,Lorton:2007,Yzelman:2011}) and allows to apply some optimization techniques, such as register blocking or cache blocking, even for higher block sizes.
\end{itemize}

Within the presented study, we consider block sizes~\eqref{eq:kl} and set $\kmax= \lmax =8$. The choice of these upper bounds stemmed from our auxiliary experiments 
which showed that space-optimal block sizes have mostly less than 64 rows/columns. Taking into account block sizes with up to 256 rows/columns should cover even the remaining corner cases.

In the summary, our optimization space is initially defined by $\schemes_6 \times \bsizes_{64}$, where $\schemes_6$ denotes a set of selected blocking storage schemes:
\begin{equation}
\schemes_6 = \bigl\{ \text{COO}, \text{CSR}, \text{bitmap}, \text{dense}, \text{min-fixed}, \text{adaptive} \bigr\}
\end{equation}
and $\bsizes_{64}$ denotes a set of selected block sizes:
\begin{equation}
\bsizes_{64} = \bigl\{ 2^k\times 2^\elll : 1\leq k,\elll \leq 8 \bigr\}.
\end{equation}

Additionally, when measuring matrix memory footprints, we need to decide how to represent information about nonzero blocks and how to represent indices. In the presented study, we assume:
\begin{enumerate}
\item nonzero blocks stored in memory in the lexicographical order;
\item explicit storage of block column index for each nonzero block;
\item storage of the number of nonzero blocks for each block row;
\item a minimum possible number of bits, i.e., $\lceil \log_2 n\rceil$ bits, to store an index related to $n$ entities (such an approach is in the literature sometimes referred to as \emph{index compression}).
\end{enumerate}

Sparse matrices are often divided into two main categories---\emph{high performance computing (HPC) matrices} and \emph{graph matrices}, the latter being binary matrices for unweighted graphs. Efficient processing of graph matrices is generally governed by special rules that are different from those being effective for HPC matrices~\cite{Ashari:2014,Buono:2016,Yang:2011} (e.g., higher matrix memory footprints in some cases lead to higher performance of computations and graph matrices are also typically not suitable for simple block processing mainly due to emergence of hypersparse blocks~\cite{Buluc:2011,Buono:2016}). Within this work, we focused mainly (but not exclusively) on HPC matrices. Particularly, for experiments, we took real matrices from the UFSMC that contained more than $10^5$ nonzero elements and that exhibited a unique structure of nonzero elements.%
\footnote{As for April, 2016.}
 This way, we obtained 563 sparse matrices arising from different application problems (see \autoref{tab:problems}) and thus having different structural (and numerical) properties; we denote these matrices by $A_1,\ldots,A_{563}$. Of these matrices, 281 were square symmetric and the remaining 282 were either rectangular or square unsymmetric.

For symmetric matrices, we always assume storage only of their single triangular parts in memory, which is a common practice. When referring to the \emph{number of nonzero elements} of a matrix, we thus generally need to distinguish between the number of \emph{all} nonzero elements and the number of elements that are assumed to be \emph{stored} in a computer memory. While measuring memory footprints of sparse matrices, we take into account the latter one.

\begin{table}[t]
\centering
\small
\begin{tabular}{lr}
\toprule
\multicolumn{1}{l}{Problem} & \multicolumn{1}{c}{Matrices} \\ \cmidrule(r){1-1} \cmidrule(l){2-2}
2D/3D                         & 36 \\
acoustics                     &  4 \\
chemical process simulation   & 25 \\
circuit simulationi           & 41 \\
computational fluid dynamics  & 47 \\
computer graphics/vision      &  8 \\
counter-example               &  2 \\
duplicate model reduction     &  5 \\
economic                      & 24 \\
eigenvalue/model reduction    &  2 \\
electromagnetics              & 11 \\
frequency-domain circuit sim. &  4 \\
\bottomrule
\end{tabular}
\hspace{6pt}
\begin{tabular}{lr}
\toprule
\multicolumn{1}{l}{Problem} & \multicolumn{1}{c}{Matrices} \\ \cmidrule(r){1-1} \cmidrule(l){2-2}
least squares                 &  7 \\
linear programming            & 51 \\
materials                     & 15 \\
model reduction               & 11 \\
optimization                  & 66 \\
power network                 & 35 \\
semiconductor device          & 16 \\
statistical/mathematical      &  1 \\
structural                    & 82 \\
theoretical/quantum chemistry & 42 \\
thermal                       & 11 \\
weighted graph                & 17 \\
\bottomrule
\end{tabular}
\caption{Counts of tested matrices falling under particular problem types (referred to as ``kinds'' in the UFSMC).}
\label{tab:problems}
\end{table}


According to the text above, a matrix memory footprint for a sparse matrix $A_k$ partitioned into uniformly-sized blocks is a function of the following parameters:
\begin{enumerate}
\item a sparse matrix itself ($A_k$),
\item a block storage scheme $\scheme\in \schemes_6$,
\item a block size $h\times w\in \bsizes_{64}$,
\item a number of bits $b$ required to store a value of a single matrix nonzero element.
\end{enumerate}
We denote this function by $\mmf(A_k,s,w\times h,b)$. We further assume storing values of matrix nonzero elements in either single or double precision IEEE floating-point format~\cite{IEEE754:2008}, which implies $b=32$ or $b=64$, respectively, in case of real matrices. We refer to such a floating-point precision as \emph{precision} only.

We say that a matrix memory footprint for a given matrix $A$ and a given precision determined by $b$ is \emph{optimal} (with respect to our work) if it equals
\begin{equation}
\min \bigl\{ \mmf(A, \scheme, h\times w, b) : \scheme\in\schemes_6, h\times w\in\bsizes_{64} \bigr\}.
\label{eq:minmmf}
\end{equation}
We call the corresponding blocking storage scheme and block size optimal as well.

Let $\schemes\subseteq \schemes_6$ and $\bsizes\subseteq \bsizes_{64}$. $\schemes\times \bsizes$ thus define a subspace of the optimization space $\schemes_6\times \bsizes_{64}$. Let
\begin{equation}
\Delta_{\schemes,\bsizes}^b(k) =
\left( \frac{\min \bigl\{ \mmf(A_k, s, h\times w, b) : s\in \schemes, h\times w \in \bsizes \bigr\}}
{ \min \bigl\{ \mmf(A_k, s, h\times w, b) : s\in \schemes_6, h\times w \in \bsizes_{64} \bigr\} }
- 1 \right) \times 100.
\end{equation}
This function expresses of how much percent is the minimal memory footprint of $A_k$ from $\schemes\times \bsizes$ higher (worse) than its optimal memory footprint. To assess the subspace $\schemes\times \bsizes$, we define the following parametrized set
\begin{equation}
\uset{\schemes}{\bsizes}{b} = \left\{ \Delta_{\schemes,\bsizes}^b(k) : 1 \leq k \leq 563 \right\}.
\end{equation}
The minimum, mean (average; $\mu$), and maximum of \uset{\schemes}{\bsizes}{b} then reflect the best, average, and worst cases, respectively, for $\schemes\times \bsizes$ across the tested matrices. If \schemes or \bsizes consists of a single element only, we omit the curly braces in the subscript of $\mathcal{U}$ for the sake of readability; e.g., we write $\uset{s}{\bsizes_{64}}{b}$ and $\uset{\schemes_6}{h\times w}{b}$ instead of $\uset{\{s\}}{\bsizes_{64}}{b}$ and $\uset{\schemes_6}{\{h\times w\}}{b}$.

\section{Results and Discussion}

First, we assessed blocking storage schemes. \autoref{tab:optsch} shows for how many tested matrices were individual schemes optimal. The adaptive scheme clearly dominates this evaluation metric; it was optimal for 464 tested matrices, which corresponds to \SI{82.4}{\percent} of their total count. Note that the min-fixed scheme was never optimal; this is due to the necessity to store additional information about the format used for blocks (if we ignored the additional 2 bits required by this scheme, it would be optimal for $58+36+5=99$ matrices). However, the numbers in \autoref{tab:optsch} reflect only best cases, i.e., matrices that were most suitable for particular schemes. To find out how much were particular schemes better than the others in average and for their worst-case (most unsuitable) matrices, we need complete statistics of \uset{s}{\bsizes_{64}}{b}; these are presented in \autoref{tab:schstats} and lead to the following observations: 
\begin{itemize}[---]
\item No fixed-format scheme minimized matrix memory footprints in comparison with the others. Bitmap was the best in average, however, it was inferior to both COO and CSR in worst cases.
\item Dense provided extremely high matrix memory footprints in average and worst cases. Due to the explicit storage of zero elements, this scheme is suitable only for kinds of matrices that contain highly dense blocks; obviously, there were only few such matrices in our tested suite (recall that the dense scheme was optimal for 5 matrices according to \autoref{tab:optsch}).
\item The lowest memory footprints were provided by the min-fixed and adaptive schemes; their numbers are considerably lower in comparison with the fixed-format schemes.
\end{itemize}
%

\begin{table}[t]
\centering
\small
\begin{tabular}{lr}
\toprule
Scheme & \multicolumn{1}{c}{Matrices} \\ \cmidrule(r){1-1} \cmidrule(l){2-2}
COO & 58 \\
CSR & 0 \\
bitmap & 36 \\
dense & 5 \\
min-fixed & 0 \\
adaptive & 464 \\
\bottomrule
\end{tabular}
\caption{Counts of tested matrices for which are blocking storage schemes optimal; the numbers are the same for both single and double precision.}
\label{tab:optsch}
\end{table}

\begin{table}[t]
\centering
\small
\begin{tabular}{l*{6}{d{-2}}}
\toprule
& \multicolumn{3}{c}{Single precision ($b=32$)} & \multicolumn{3}{c}{Double precision ($b=64$)} \\
	\cmidrule(l){2-4} \cmidrule(l){5-7} 
Scheme ($s$)
& \multicolumn{1}{c}{Minimum} & \multicolumn{1}{c}{Average} & \multicolumn{1}{c}{Maximum}
& \multicolumn{1}{c}{Minimum} & \multicolumn{1}{c}{Average} & \multicolumn{1}{c}{Maximum} \\
	\cmidrule(r){1-1} \cmidrule(l){2-4} \cmidrule(l){5-7}
COO & 0.00 & 4.78 & 15.27 & 0.00 & 2.52 & 7.67 \\
CSR & 0.73 & 6.84 & 19.13 & 0.41 & 3.74 & 11.05 \\
bitmap & 0.00 & 3.13 & 22.01 & 0.00 & 1.75 & 12.38 \\
dense & 0.00 & 84.61 & 217.04 & 0.00 & 92.40 & 249.02 \\
min-fixed & 0.00 & 1.19 & 5.41 & 0.00 & 0.64 & 2.94 \\
adaptive & 0.00 & 0.10 & 2.24 & 0.00 & 0.05 & 1.30 \\
\bottomrule
\end{tabular}
\caption{Minimum, average and maximum values of $\uset{s}{\bsizes_{64}}{b}$ (in percents).}
\label{tab:schstats}
\end{table}

\begin{figure}[t]
\centering
\includegraphics[width=\textwidth]{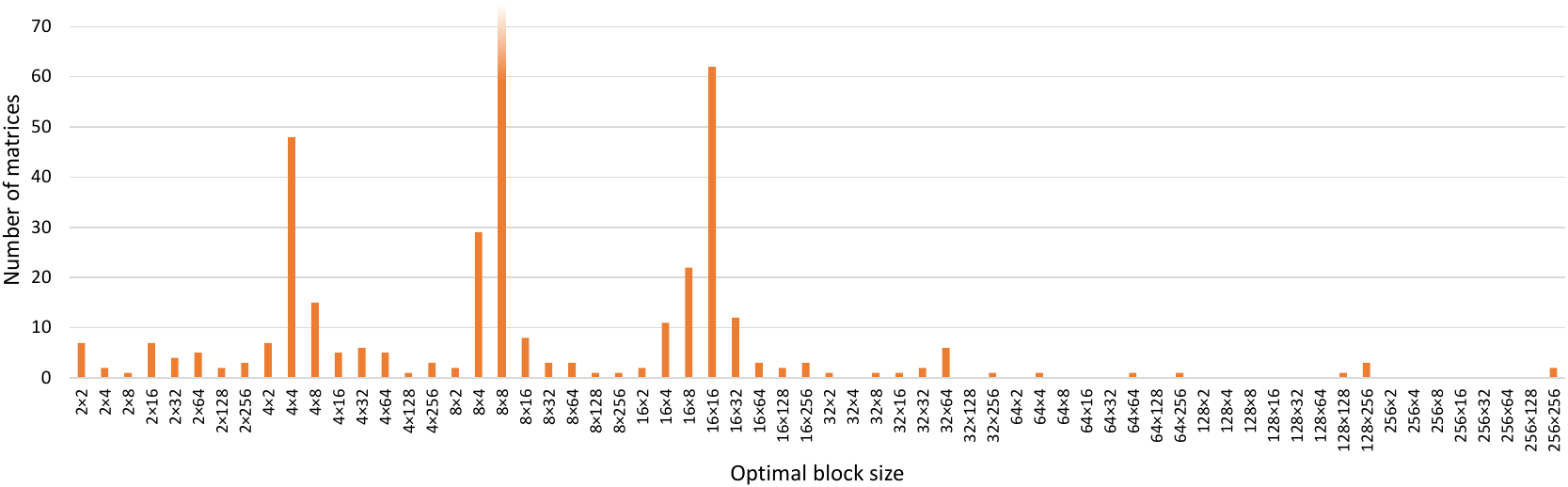}
\caption{Numbers of tested matrices for which are block sizes optimal, measured for double precision; block size $8\times 8$ was optimal for  257 matrices.}
\label{fig:optbs}
\end{figure}

Similarly as blocking storage schemes, we assessed block sizes. \autoref{fig:optbs} shows for how many tested matrices were individual block sizes optimal in case of double precision measurements; for single precision, the results differed only for 2 matrices. 
We may observe that some block sizes were especially favourable. The $8\times 8$ block size was optimal for 257 matrices, which corresponds to \SI{45.6}{\percent} of their total count. Together with $4\times 4$ and $16\times 16$, these 3 block sizes were optimal for \SI{65.2}{\percent} of tested matrices. However, again, the numbers from \autoref{fig:optbs} reflect only best cases. To find out how much were particular block sizes better than the others in average and for their worst-cases matrices, we present the average and maximum values of \uset{\schemes_6}{h\times w}{b} in \autoref{tab:bsstats:single} and \autoref{tab:bsstats:double} for single and double precision, respectively. According to these results, some blocks sizes---especially $8\times 8$---provided alone average matrix memory footprints close to their optimal values. However, there was not a single block size that would yield the same outcome for all the tested matrices; the maxima were for all the block sizes relatively high.


\addtolength{\tabcolsep}{-3pt}
\begin{table}[t]
\centering
\small
\begin{tabular}[t]{rx{-3}d{-2}d{-2}}
\toprule
\multicolumn{1}{c}{Rank} & \multicolumn{1}{c}{$h\times w$} &
	\multicolumn{1}{c}{Avg.} & \multicolumn{1}{c}{Max.} \\
	\cmidrule(r){1-2} \cmidrule(l){3-4} 
1 & 8x8 & 1.23 & 18.36 \\
2 & 8x16 & 2.14 & 19.35 \\
3 & 16x8 & 2.26 & 21.41 \\
4 & 4x8 & 2.32 & 17.31 \\
5 & 8x4 & 2.38 & 19.52 \\
6 & 16x16 & 2.56 & 21.82 \\
7 & 4x4 & 2.92 & 21.94 \\
8 & 4x16 & 2.99 & 16.51 \\
9 & 16x4 & 3.23 & 20.44 \\
10 & 8x32 & 3.65 & 21.26 \\
\bottomrule
\end{tabular}
\hspace{6pt}
\begin{tabular}[t]{rx{-3}d{-2}d{-2}}
\toprule
\multicolumn{1}{c}{Rank} & \multicolumn{1}{c}{$h\times w$} &
	\multicolumn{1}{c}{Avg.} & \multicolumn{1}{c}{Max.} \\
	\cmidrule(r){1-2} \cmidrule(l){3-4} 
11 & 16x32 & 4.03 & 23.75 \\
12 & 32x8 & 4.13 & 23.97 \\
13 & 4x32 & 4.36 & 18.71 \\
14 & 32x16 & 4.53 & 24.45 \\
15 & 32x4 & 4.87 & 23.60 \\
16 & 32x32 & 5.20 & 26.50 \\
17 & 2x8 & 5.59 & 21.15 \\
18 & 8x64 & 5.61 & 23.57 \\
19 & 8x2 & 5.66 & 26.39 \\
20 & 2x16 & 5.84 & 22.84 \\
\bottomrule
\end{tabular}
\hspace{6pt}
\begin{tabular}[t]{rx{-3}d{-2}d{-2}}
\toprule
\multicolumn{1}{c}{Rank} & \multicolumn{1}{c}{$h\times w$} &
	\multicolumn{1}{c}{Avg.} & \multicolumn{1}{c}{Max.} \\
	\cmidrule(r){1-2} \cmidrule(l){3-4} 
21 & 16x64 & 5.89 & 26.15 \\
22 & 4x2 & 6.06 & 28.77 \\
23 & 2x4 & 6.15 & 23.07 \\
24 & 16x2 & 6.25 & 29.98 \\
25 & 4x64 & 6.26 & 21.53 \\
26 & 64x8 & 6.56 & 25.83 \\
\multicolumn{1}{r}{\dots} & \multicolumn{1}{c}{\dots} & \multicolumn{1}{r}{\dots} & \multicolumn{1}{r}{\dots} \\
62 & 256x2 & 14.44 & 37.33 \\
63 & 256x128 & 14.61 & 38.32 \\
64 & 256x256 & 14.65 & 35.42 \\
\bottomrule
\end{tabular}
\caption{Average and maximum values of \uset{\schemes_6}{h\times w}{32} (in percents), sorted by average.}.
\label{tab:bsstats:single}
\end{table}
\addtolength{\tabcolsep}{3pt}

\addtolength{\tabcolsep}{-3pt}
\begin{table}[t]
\centering
\small
\begin{tabular}[t]{rx{-3}d{-2}d{-2}}
\toprule
\multicolumn{1}{c}{Rank} & \multicolumn{1}{c}{$h\times w$} &
	\multicolumn{1}{c}{Avg.} & \multicolumn{1}{c}{Max.} \\
	\cmidrule(r){1-2} \cmidrule(l){3-4} 
1 & 8x8 & 0.69 & 11.07 \\
2 & 8x16 & 1.18 & 11.67 \\
3 & 16x8 & 1.25 & 12.91 \\
4 & 4x8 & 1.30 & 9.74 \\
5 & 8x4 & 1.33 & 10.98 \\
6 & 16x16 & 1.40 & 13.16 \\
7 & 4x4 & 1.63 & 12.34 \\
8 & 4x16 & 1.66 & 9.96 \\
9 & 16x4 & 1.79 & 12.32 \\
10 & 8x32 & 1.99 & 11.97 \\
\bottomrule
\end{tabular}
\hspace{6pt}
\begin{tabular}[t]{rx{-3}d{-2}d{-2}}
\toprule
\multicolumn{1}{c}{Rank} & \multicolumn{1}{c}{$h\times w$} &
	\multicolumn{1}{c}{Avg.} & \multicolumn{1}{c}{Max.} \\
	\cmidrule(r){1-2} \cmidrule(l){3-4} 
11 & 16x32 & 2.19 & 12.84 \\
12 & 32x8 & 2.26 & 14.45 \\
13 & 4x32 & 2.40 & 10.56 \\
14 & 32x16 & 2.47 & 14.04 \\
15 & 32x4 & 2.68 & 14.23 \\
16 & 32x32 & 2.82 & 14.18 \\
17 & 8x64 & 3.05 & 12.62 \\
18 & 2x8 & 3.11 & 12.08 \\
19 & 8x2 & 3.14 & 14.02 \\
20 & 16x64 & 3.19 & 14.00 \\
\bottomrule
\end{tabular}
\hspace{6pt}
\begin{tabular}[t]{rx{-3}d{-2}d{-2}}
\toprule
\multicolumn{1}{c}{Rank} & \multicolumn{1}{c}{$h\times w$} &
	\multicolumn{1}{c}{Avg.} & \multicolumn{1}{c}{Max.} \\
	\cmidrule(r){1-2} \cmidrule(l){3-4} 
21 & 2x16 & 3.25 & 13.04 \\
22 & 4x2 & 3.34 & 15.74 \\
23 & 2x4 & 3.40 & 12.84 \\
24 & 4x64 & 3.42 & 11.38 \\
25 & 16x2 & 3.47 & 15.93 \\
26 & 64x8 & 3.57 & 15.30 \\
\multicolumn{1}{r}{\dots} & \multicolumn{1}{c}{\dots} & \multicolumn{1}{r}{\dots} & \multicolumn{1}{r}{\dots} \\
62 & 256x2 & 7.88 & 21.59 \\
63 & 256x128 & 7.92 & 19.56 \\
64 & 256x256 & 7.93 & 18.96 \\
\bottomrule
\end{tabular}
\caption{Average and maximum values of \uset{\schemes_6}{h\times w}{64} (in percents), sorted by average.}.
\label{tab:bsstats:double}
\end{table}
\addtolength{\tabcolsep}{3pt}

\begin{figure}[t]
\centering
\includegraphics[width=\textwidth]{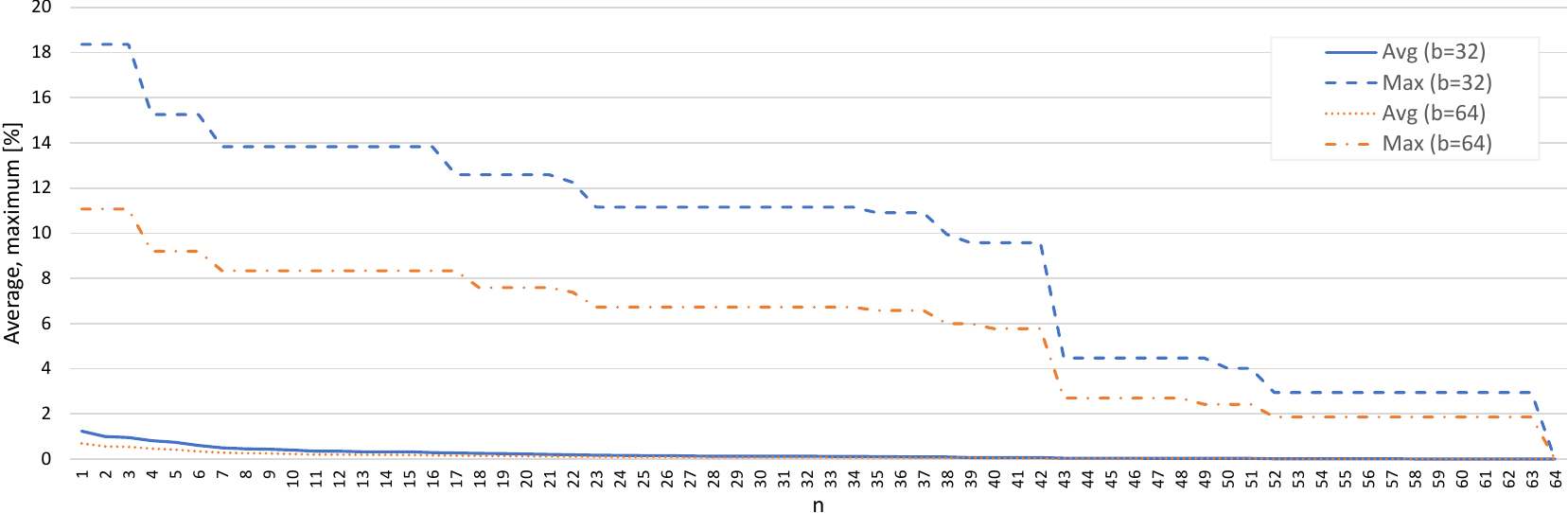}
\caption{Average and maximum values \uset{\schemes_6}{\cset{n}{b}}{b} (in percents) as a funciton of $n$.}
\label{fig:cnstats}
\end{figure}

Let us remind that one of our goals is a possible reduction of the number of block sizes in the optimization test space. 
The question thus is whether there is some subset $\bsizes\subset \bsizes_{64}$ that would, at the same time:
\begin{enumerate}
\item significantly reduce the number of block sizes ($\lvert \bsizes \rvert$),
\item provide matrix memory footprints close to their optimal values for most of the tested matrices (average of $\uset{\schemes_6}{B}{b}$ close to zero),
\item provide low matrix memory footprints for all the tested matrices (low maximum of $\uset{\schemes_6}{B}{b}$).
\end{enumerate}
Natural candidates for such a subset would be the first $n$ block sizes from \autoref{tab:bsstats:single} and \autoref{tab:bsstats:double}; let us denote them by \cset{n}{64} and \cset{n}{32}, respectively. \autoref{fig:cnstats} evaluates these subsets as a function of $n$. We may notice that 
\begin{align}
\cset{9}{64} =\cset{9}{32} & = \bigl\{\, h\times w : h,w \in \{4,8,16\} \bigr\}, \\ 
\cset{16}{64} =\cset{16}{32} & = \bigl\{\, h\times w : h,w \in \{4,8,16,32\} \bigr\}; 
\end{align} 
seemingly, block sizes from these subsets are especially suitable for sparse matrices in general.


Despite that, neither these first 9 nor 16 block sizes reduced the maximal matrix memory footprints too much according to \autoref{fig:cnstats}. However, we may observe that there are some block sizes where these maxima 
significantly dropped. Based on the analysis of the statistics of \uset{\schemes_6}{\cset{n}{b}}{b}, we propose the following \emph{reduced sets of block sizes}:
\begin{align}
\bsizes_8    &= \bigl\{ 2^k\times 2^k : 1\leq k\leq 8 \bigr\}, \\
\bsizes_{14} &= \bsizes_8 \cup \bigl\{2^k\times 2^\elll : 2\leq k,\elll \leq 4\bigr\}, \\
\bsizes_{20} &= \bsizes_8 \cup \bigl\{2^k\times 2^\elll : 2\leq k,\elll \leq 5\bigr\}.
\end{align}
$\bsizes_8$ thus consists of all square block sizes from $\bsizes_{64}$. $\bsizes_{14}$ and $\bsizes_{20}$ equal $\bsizes_8$ plus rectangular block sizes from \cset{9}{32} (\cset{9}{64}) and \cset{16}{32} (\cset{16}{64}), respectively.

\autoref{tab:schstats} revealed that to minimize memory footprints of (all) the tested matrices, we had to use either the min-fixed or the adaptive blocking storage scheme. To reduce the block preprocessing overhead, we now proposed several reduced sets of block sizes. Let us now assess these options together. We measured the statistics of \uset{s}{\bsizes_j}{b} for all the combinations of $s\in\{\text{min-fixed},\text{adaptive}\}$ and $j\in \{64,20,14,8\}$; the results are presented in \autoref{tab:reduced}. The average matrix memory footprints were in all cases close to their optimal values. Moreover, the reduced sets $\bsizes_j$ required much less block sizes than \cset{n}{b} to achieve the same maxima. For instance:
\begin{enumerate}
\item $\bsizes_{14}$ in combination with the min-fixed scheme required only 14 block sizes to achieve the same maxima as \cset{43}{b} in combination with all the schemes. This would effectively reduce the number of block sizes in the optimization space by a factor of about 3, which would proportionally reduce the preprocessing overhead in practice.

\item $\bsizes_{20}$ in combination with the adaptive scheme required only 20 block sizes to achieve the same maxima \cset{50}{b} in combination with all the schemes. This would effectively reduce the number of block sizes by a factor of $2.5$.
\end{enumerate}

\begin{table}[t]
\centering
\small
\begin{subtable}{\textwidth}
\centering
\caption{Single precision ($b=32$)}
\begin{tabular}{c*{4}{d{-3}}}
\toprule
& \multicolumn{2}{c}{$s=\text{min-fixed}$} & \multicolumn{2}{c}{$s=\text{adaptive}$} \\
	\cmidrule(l){2-3} \cmidrule(l){4-5} 
Block sizes
& \multicolumn{1}{c}{Average} & \multicolumn{1}{c}{Maximum}
& \multicolumn{1}{c}{Average} & \multicolumn{1}{c}{Maximum} \\
	\cmidrule(r){1-1} \cmidrule(l){2-3} \cmidrule(l){4-5}
$\bsizes_{64}$ & 1.19 &  5.41 & 0.10 &  2.24 \\
$\bsizes_{20}$ & 1.32 &  6.23 & 0.22 &  4.21 \\
$\bsizes_{14}$ & 1.35 &  6.89 & 0.28 &  6.81 \\
$\bsizes_8$    & 1.51 & 10.06 & 0.51 & 11.07 \\
\bottomrule
\end{tabular}
\end{subtable}
\\ \vspace{1em}
\begin{subtable}{\textwidth}
\centering
\caption{Double precision ($b=64$)}
\begin{tabular}{c*{4}{d{-3}}}
\toprule
& \multicolumn{2}{c}{$s=\text{min-fixed}$} & \multicolumn{2}{c}{$s=\text{adaptive}$} \\
	\cmidrule(l){2-3} \cmidrule(l){4-5} 
Block sizes
& \multicolumn{1}{c}{Average} & \multicolumn{1}{c}{Maximum}
& \multicolumn{1}{c}{Average} & \multicolumn{1}{c}{Maximum} \\
	\cmidrule(r){1-1} \cmidrule(l){2-3} \cmidrule(l){4-5}
$\bsizes_{64}$ & 0.64 & 2.94 & 0.05 & 1.30 \\
$\bsizes_{20}$ & 0.71 & 3.52 & 0.12 & 2.37 \\
$\bsizes_{14}$ & 0.73 & 3.77 & 0.16 & 3.83 \\
$\bsizes_8$    & 0.81 & 5.34 & 0.28 & 5.88 \\
\bottomrule
\end{tabular}
\end{subtable}
\caption{Average and maximum values of \uset{s}{\bsizes_j}{b} (in percents) for $j\in\{64,20,14,8\}$.}
\label{tab:reduced}
\end{table}

\subsection{Consistency}

Up to now, we have presented measurements conducted for all 563 tested matrices. To asses their ``representativeness'', we measured the consistency of memory footprints statistics across randomly selected subsets of these matrices. Such an experiment should reveal how our measurements are sensitive to the set of input matrices, which should suggest to which extent we can generalize the outcomes of the study.

Let \aset{n}{i} denote an $i$th set of $n$ randomly selected tested matrices; different $i$ thus allows us to distinguish different random selections. Let \kset{n}{i} denote a set of matrix indices from \aset{n}{i}, thus $\aset{n}{i}=\bigl\{ A_k : k \in \kset{n}{i} \bigr\}$. Let 
\begin{equation}
\vset{s}{\bsizes_j}{b}{n}{i} = \left\{ \Delta_{s,\bsizes_j}^b(k) : k \in \kset{n}{i} \right\}.
\end{equation}
\vset{s}{\bsizes_j}{b}{n}{i} thus expresses of how much percents are memory footprints of matrices from \aset{n}{i}---measured for scheme $s$, a set of block sizes $\bsizes_j$, and a precision given by $b$---higher than their optimal memory footprints. Similarly as before, we were interested in average and maximum values of \vset{s}{\bsizes_j}{b}{n}{i}; let them denote by $\avg \vset{s}{\bsizes_j}{b}{n}{i}$ and $\max \vset{s}{\bsizes_j}{b}{n}{i}$, respectively. To assess the consistency introduced above, we measured standard deviations of these metrics for 50 sets of 200 randomly selected tested matrices, i.e., standard deviations of the following sets:
\begin{equation}
\Bigl\{ \avg \vset{s}{\bsizes_j}{b}{200}{i} : 1 \leq i \leq 50 \Bigr\}
\quad \text{and} \quad
\Bigl\{ \max \vset{s}{\bsizes_j}{b}{200}{i} : 1 \leq i \leq 50 \Bigr\}.
\end{equation}
The results obtained for the min-fixed and adaptive schemes, sets of blocks sizes $\bsizes_{64},\bsizes_{20},\bsizes_{14},\bsizes_8$, and both precisions are shown in \autoref{tab:consistency}.

The measured standard deviations are of 1 to 2 orders of magnitude lower than the corresponding numbers from \autoref{tab:reduced}. By normalizing the standard deviations by these numbers, we found out that the standard deviations ranged from $5.16$ to $9.28$ percents for the min-fixed scheme and from $10.30$ to $21.33$ percents for the adaptive scheme. Seemingly, the min-fixed scheme provides more consistent relative memory footprints of matrices with respect to their optimal values, while the adaptive scheme is more sensitive to the selection of matrices as for this evaluation metric. Note, however, that the measured standard deviations were according to \autoref{tab:consistency} in all cases relatively small with the maximum value $1.41$; recall that these numbers are relative differences in percents between optimal matrix memory footprints and those measured for particular tested configurations. Especially, the standard deviations for average metrics are practically negligible, which manifests high level of representativeness of the tested matrices.

\begin{table}[t]
\centering
\small
\begin{subtable}{\textwidth}
\centering
\caption{Single precision ($b=32$)}
\begin{tabular}{c*{4}{d{-3}}}
\toprule
& \multicolumn{2}{c}{$s=\text{min-fixed}$} & \multicolumn{2}{c}{$s=\text{adaptive}$} \\
	\cmidrule(l){2-3} \cmidrule(l){4-5} 
Block sizes
& \multicolumn{1}{c}{Average} & \multicolumn{1}{c}{Maximum}
& \multicolumn{1}{c}{Average} & \multicolumn{1}{c}{Maximum} \\
	\cmidrule(r){1-1} \cmidrule(l){2-3} \cmidrule(l){4-5}
$\bsizes_{64}$ & 0.06 &  0.33 & 0.02 &  0.24 \\
$\bsizes_{20}$ & 0.07 &  0.37 & 0.03 &  0.64 \\
$\bsizes_{14}$ & 0.08 &  0.42 & 0.04 &  1.41 \\
$\bsizes_8$    & 0.09 &  0.86 & 0.07 &  1.40 \\
\bottomrule
\end{tabular}
\end{subtable}
\\ \vspace{1em}
\begin{subtable}{\textwidth}
\centering
\caption{Double precision ($b=64$)}
\begin{tabular}{c*{4}{d{-3}}}
\toprule
& \multicolumn{2}{c}{$s=\text{min-fixed}$} & \multicolumn{2}{c}{$s=\text{adaptive}$} \\
	\cmidrule(l){2-3} \cmidrule(l){4-5} 
Block sizes
& \multicolumn{1}{c}{Average} & \multicolumn{1}{c}{Maximum}
& \multicolumn{1}{c}{Average} & \multicolumn{1}{c}{Maximum} \\
	\cmidrule(r){1-1} \cmidrule(l){2-3} \cmidrule(l){4-5}
$\bsizes_{64}$ & 0.03 & 0.19 & 0.01 & 0.14 \\
$\bsizes_{20}$ & 0.04 & 0.23 & 0.02 & 0.29 \\
$\bsizes_{14}$ & 0.04 & 0.22 & 0.02 & 0.82 \\
$\bsizes_8$    & 0.05 & 0.50 & 0.04 & 0.61 \\
\bottomrule
\end{tabular}
\end{subtable}
\caption{Standard deviations of $\avg \vset{s}{\bsizes_j}{b}{200}{i}$ and $\max \vset{s}{\bsizes_j}{b}{200}{i}$ (in percents) for $1 \leq i \leq 50$.}
\label{tab:consistency}
\end{table}

\subsection{Blocking Storage Schemes Without CSR}

We have defined the min-fixed and adaptive blocking storage schemes such that the format used for storing blocks is selected---from COO, CSR, bitmap, and dense---either for all blocks collectively or for each block separately; the corresponding results were presented by \autoref{tab:reduced}. However, we were also interested in how these results would change if we modified the min-fixed and adaptive schemes by excluding individual formats. We carried out such measurements and their results revealed that:
\begin{enumerate}
\item without the COO or bitmap format, the memory footprints of matrices grew significantly;
\item without the CSR or dense formats, the memory footprints of matrices grew negligibly;
\item without both the CSR and dense formats, the memory footprints of matrices grew negligibly as well.
\end{enumerate}

The question therefore is whether the CSR and dense formats are at all useful for storing blocks. Based on our knowledge and experience, we would not suggest to exclude the dense format. Though this format is optimal in rare cases only, it is likely the most efficient format for matrix computations. For example, multiplication of a block stored in the dense format with a corresponding vector part can be performed by invoking a relevant operation from some dense linear algebra library, such as BLAS~\cite{Dongarra:2002}. In practice, every HPC system provides at least one optimized implementation of such a library that is highly-tuned for a given hardware architecture (e.g., ATLAS, BLIS, Cray LibSci, IBM ESSL, Intel MKL, OpenBLAS, etc.).

On the contrary, CSR does not provide the same benefits as the dense format, especially when it is implemented together with index compression. Moreover, CSR is the only considered format that prescribes a fixed order of nonzero elements; consequently, it does not allow to store them in an order that might be computationally more efficient, such as the Z-Morton order. One therefore might consider excluding CSR from the min-fixed and adaptive schemes to simplify related algorithms and their implementations. We call such modified schemes \emph{min-fixed-w/o-CSR} and \emph{adaptive-w/o-CSR} and present the results for them in \autoref{tab:reducedwocsr}. Obviously, the numbers are either the same or only slightly higher than those measured for the original min-fixed and adaptive schemes; see \autoref{tab:reduced}.

\begin{table}[t]
\centering
\small
\begin{subtable}{\textwidth}
\centering
\caption{Single precision ($b=32$)}
\begin{tabular}{c*{4}{d{-3}}}
\toprule
& \multicolumn{2}{c}{$s=\text{min-fixed-w/o-CSR}$} & \multicolumn{2}{c}{$s=\text{adaptive-w/o-CSR}$} \\
	\cmidrule(l){2-3} \cmidrule(l){4-5} 
Block sizes
& \multicolumn{1}{c}{Average} & \multicolumn{1}{c}{Maximum}
& \multicolumn{1}{c}{Average} & \multicolumn{1}{c}{Maximum} \\
	\cmidrule(r){1-1} \cmidrule(l){2-3} \cmidrule(l){4-5}
$\bsizes_{64}$ & 1.20 &  5.55 & 0.15 &  2.24 \\
$\bsizes_{20}$ & 1.32 &  6.44 & 0.26 &  4.22 \\
$\bsizes_{14}$ & 1.35 &  6.89 & 0.31 &  6.82 \\
$\bsizes_8$    & 1.51 & 10.06 & 0.54 & 11.07 \\
\bottomrule
\end{tabular}
\end{subtable} 
\\ \vspace{1em}
\begin{subtable}{\textwidth}
\centering
\caption{Double precision ($b=64$)}
\begin{tabular}{c*{4}{d{-3}}}
\toprule
& \multicolumn{2}{c}{$s=\text{min-fixed-w/o-CSR}$} & \multicolumn{2}{c}{$s=\text{adaptive-w/o-CSR}$} \\
	\cmidrule(l){2-3} \cmidrule(l){4-5} 
Block sizes
& \multicolumn{1}{c}{Average} & \multicolumn{1}{c}{Maximum}
& \multicolumn{1}{c}{Average} & \multicolumn{1}{c}{Maximum} \\
	\cmidrule(r){1-1} \cmidrule(l){2-3} \cmidrule(l){4-5}
$\bsizes_{64}$ & 0.65 & 3.01 & 0.09 & 1.30 \\
$\bsizes_{20}$ & 0.71 & 3.52 & 0.15 & 2.37 \\
$\bsizes_{14}$ & 0.73 & 3.77 & 0.18 & 3.84 \\
$\bsizes_8$    & 0.82 & 5.34 & 0.30 & 5.88 \\
\bottomrule
\end{tabular}
\end{subtable}
\caption{Average and maximum values of \uset{s}{\bsizes_k}{b} (in percents) for $k\in\{64,20,14,8\}$ with excluded CSR.}
\label{tab:reducedwocsr}
\end{table}

\subsection{Memory Savings Against CSR32}
\label{sec:savings}

Likely the most widely-used storage format for sparse matrices in practice is CSR, which is supported by vast majority of software tools and libraries that work with sparse matrices. To distinguish between CSR used for blocks of partitioned matrices and CSR used for whole (not-partitioned) matrices, we call the latter CSR32, since it is typically implemented with 32-bit indices. Researchers frequently demonstrate the superiority of their algorithms and data structures (formats) by comparison with CSR32, which have become de facto an etalon in sparse-matrix research.

Comparison of memory footprints of sparse matrices partitioned into blocks and the same matrices stored in CSR32 allows us to assess our blocking approach. Let $\mmfcsr(A, b)$ denote a memory footprint of a matrix $A$ stored in memory in CSR32 with respect to a precision given by $b$. The function
\begin{equation}
\Lambda^b(k) =
\left(1 - \frac{ \min \bigl\{ \mmf(A_k, s, h\times w, b) : s\in \schemes_6, h\times w \in \bsizes_{64} \bigr\} }
{\mmfcsr(A_k, b)}
\right) \times 100
\end{equation}
then expresses how much memory in percents we would save if we stored the tested matrix $A_k$ in its optimal blocking configuration instead of in CSR32. Me measured these memory savings for all the tested matrices and processed them statistically; the results are presented by \autoref{tab:savings}. The obtained numbers speaks strongly in favour of partitioning of sparse matrices in general. Even in worst cases, our blocking approach reduced the memory footprints of matrices of \SI{25.46}{\percent} and \SI{17.08}{\percent} for single and double precision, respectively. In average, the savings were \SI{42.29}{\percent} and \SI{28.67}{\percent}, which significantly reduces the amount of data that needs to be transferred between memory and processors during computations.

\begin{table}[t]
\centering
\small
\begin{tabular}{l*{2}{d{-3}}}
\toprule
Statistics & \multicolumn{1}{c}{Single precision} & \multicolumn{1}{c}{Double precision} \\
	\cmidrule(r){1-1} \cmidrule(l){2-3}
Minimum & 25.46 & 17.08 \\
Average & 42.29 & 28.67 \\
Maximum & 50.21 & 35.86 \\
\bottomrule
\end{tabular}
\caption{Statistics of $\Lambda^b(k)$, i.e., memory savings of optimal blocking configurations against CSR32 in percents, across the tested matrices.}
\label{tab:savings}
\end{table}

\autoref{tab:savings} shows the statistics of memory savings across all the tested matrices. However, we also wanted to find out which matrices were especially suitable/unsuitable for partitioning in general. For this reason, we measured the memory saving against CSR32 also as a function the following criteria, which are commonly used to distinguish/quantify different types of sparse matrices:
\begin{enumerate}
\item application problem type,
\item relative count of matrix nonzero elements (their density),
\item uniformity of the distribution of matrix nonzero elements across its rows.
\end{enumerate}
The application problem types were introduced by \autoref{tab:problems}. As for the second criterion, we define the \emph{density} of nonzero elements for an $m\times n$ matrix $A$ with \nnz nonzero elements in percents as $\rho(A) = \nnz / (m \times n) \times 100$. Its values thus ranges from $0$ for an empty matrix to $100$ to a fully dense matrix.

Let $\rnnz(i)$ denote a number of nonzero elements of $i$th row of $A$; $\rnnz(i)$ thus ranges from 0 for empty rows to $n$ for fully dense rows. To allow a collective evaluation of matrices with different row lengths, we transform $\rnnz(i)$ into relative counts in percents as follows: $\prnnz(i) = \rnnz(i) / n \times 100$. The standard deviation of $\prnnz(i)$ for $i=1,\ldots,m$ then represents an inverse measure of the above introduced third criterion for $A$. Zero standard deviation of $\prnnz(i)$ then implies a matrix whose all rows have exactly the same number of nonzero elements.

Recall that in \autoref{sec:meth} we defined two kinds of the numbers of nonzero elements, counting either all of them or just those stored in a computer memory
(for unsymmetric matrices, these numbers would be equal). Accordingly, we can quantify the above introduced second and third matrix criteria in two ways; we further show results for both of them.

\begin{figure}[t]
\centering
\includegraphics[width=13cm]{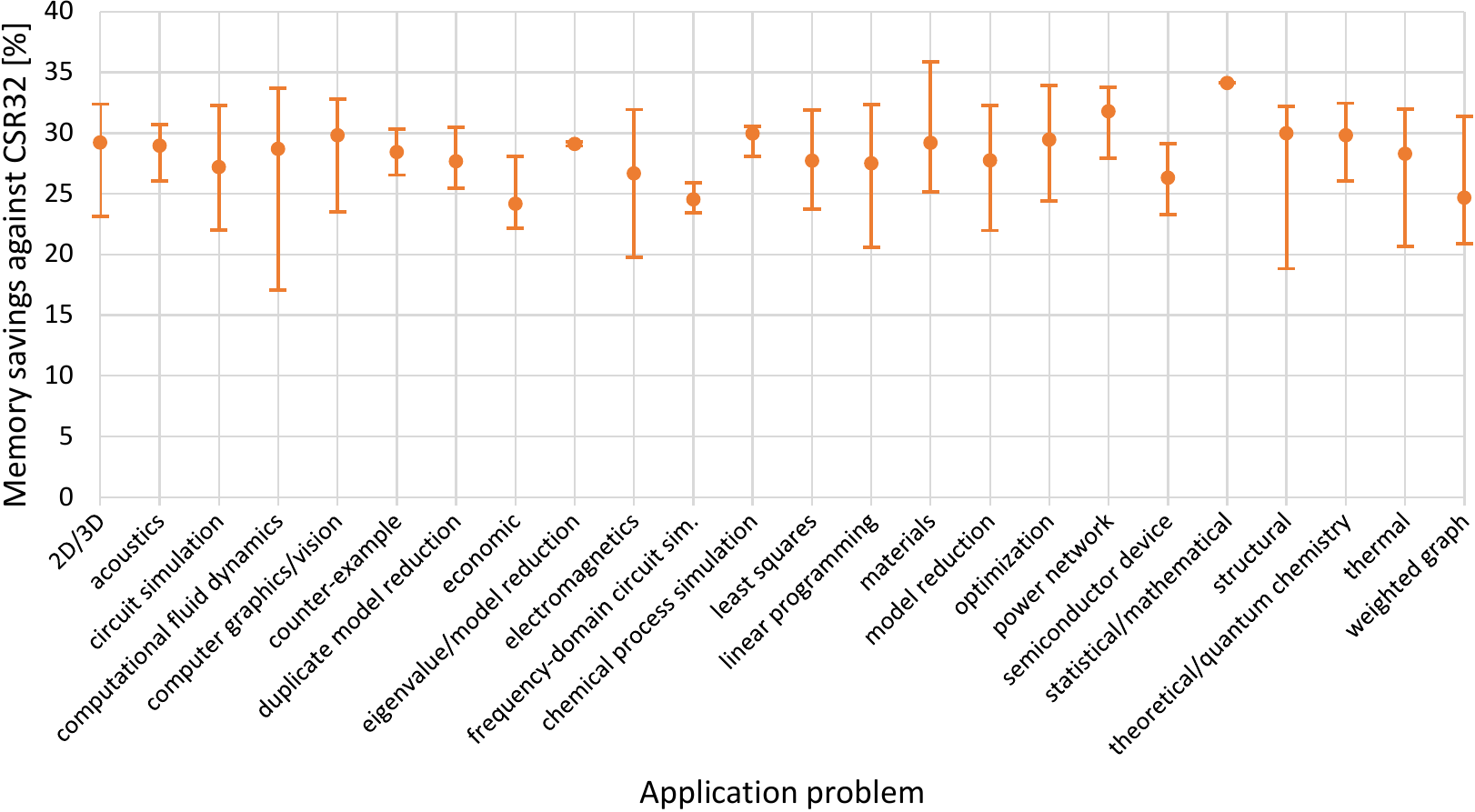}
\caption{Statistics of relative memory savings against CSR32 in percents across the tested matrices grouped by individual problem types, measured for double precision. Circles represent average values, the extents from minimal to maximal values are indicated by bars.}
\label{fig:savings_problem}
\end{figure}

The measurements for the first criterion and double precision are presented by \autoref{fig:savings_problem}; the results for single precision are practically the same, just scaled accordingly. We need to be careful when making general conclusions based on these results, since for some problem types, our tested suite of matrices contain only few representatives. However, we may observe that the memory savings against CSR32 were relatively consistent across problem types; there was no problem type that would provide much better or much worse savings than the others, including even the graph matrices. 

\begin{figure}[t]
\centering
\begin{subfigure}{.47\textwidth}
	\centering
	\includegraphics[width=\textwidth]{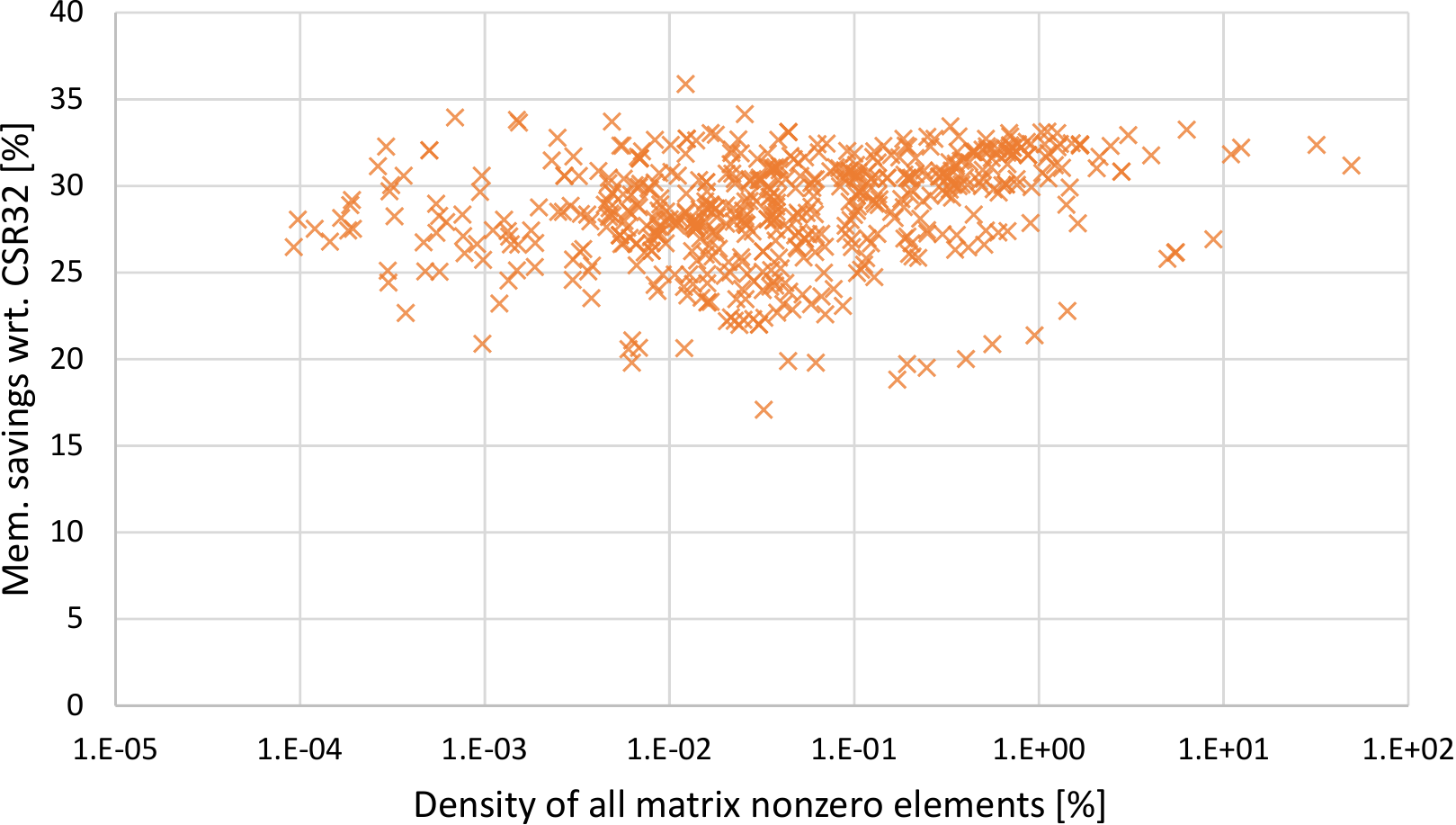}
\end{subfigure}
\hspace{1em}
\begin{subfigure}{.47\textwidth}
	\centering
	\includegraphics[width=\textwidth]{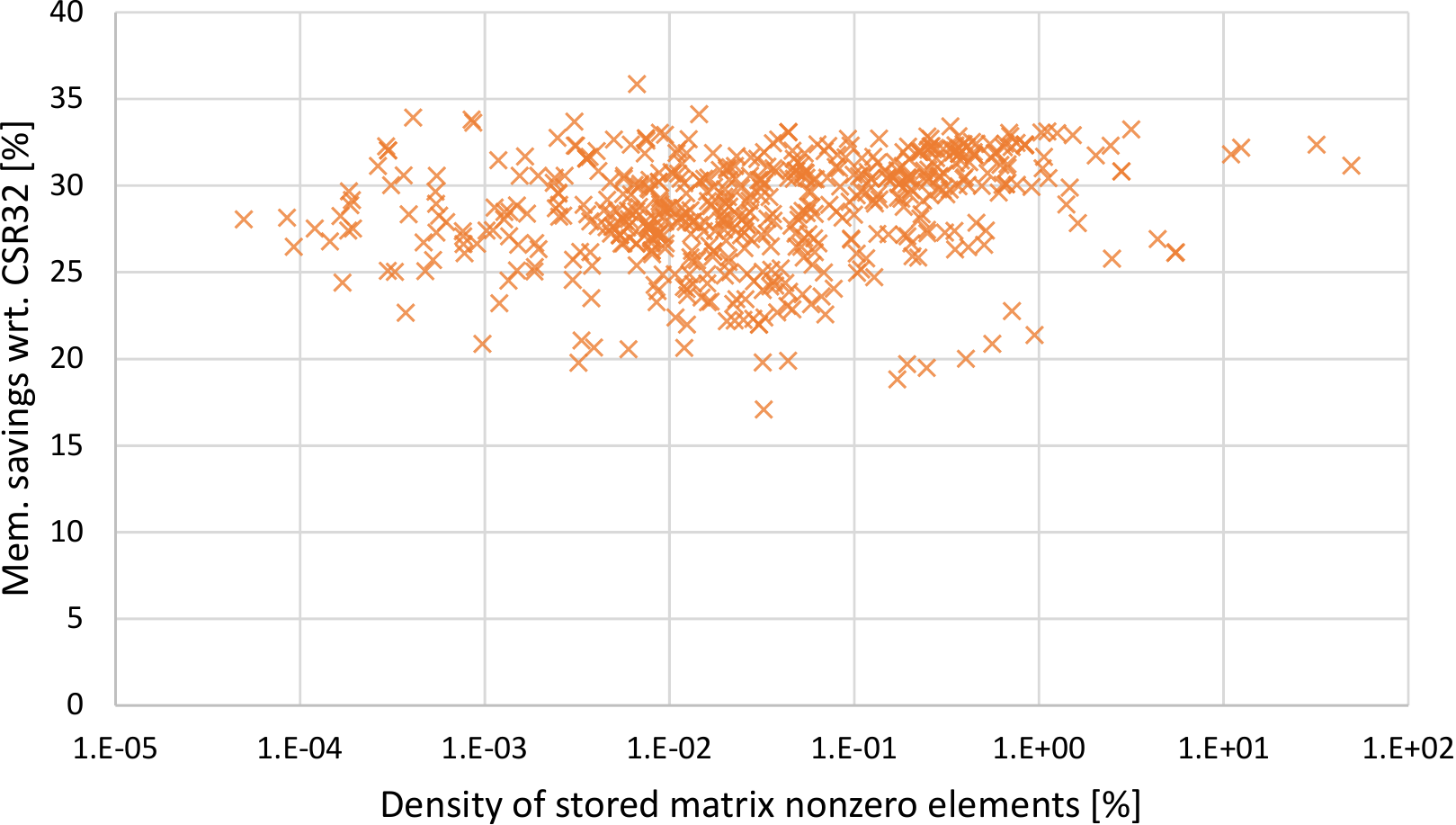}
\end{subfigure}
\\ \vspace{1em}
\begin{subfigure}{.47\textwidth}
	\centering
	\includegraphics[width=\textwidth]{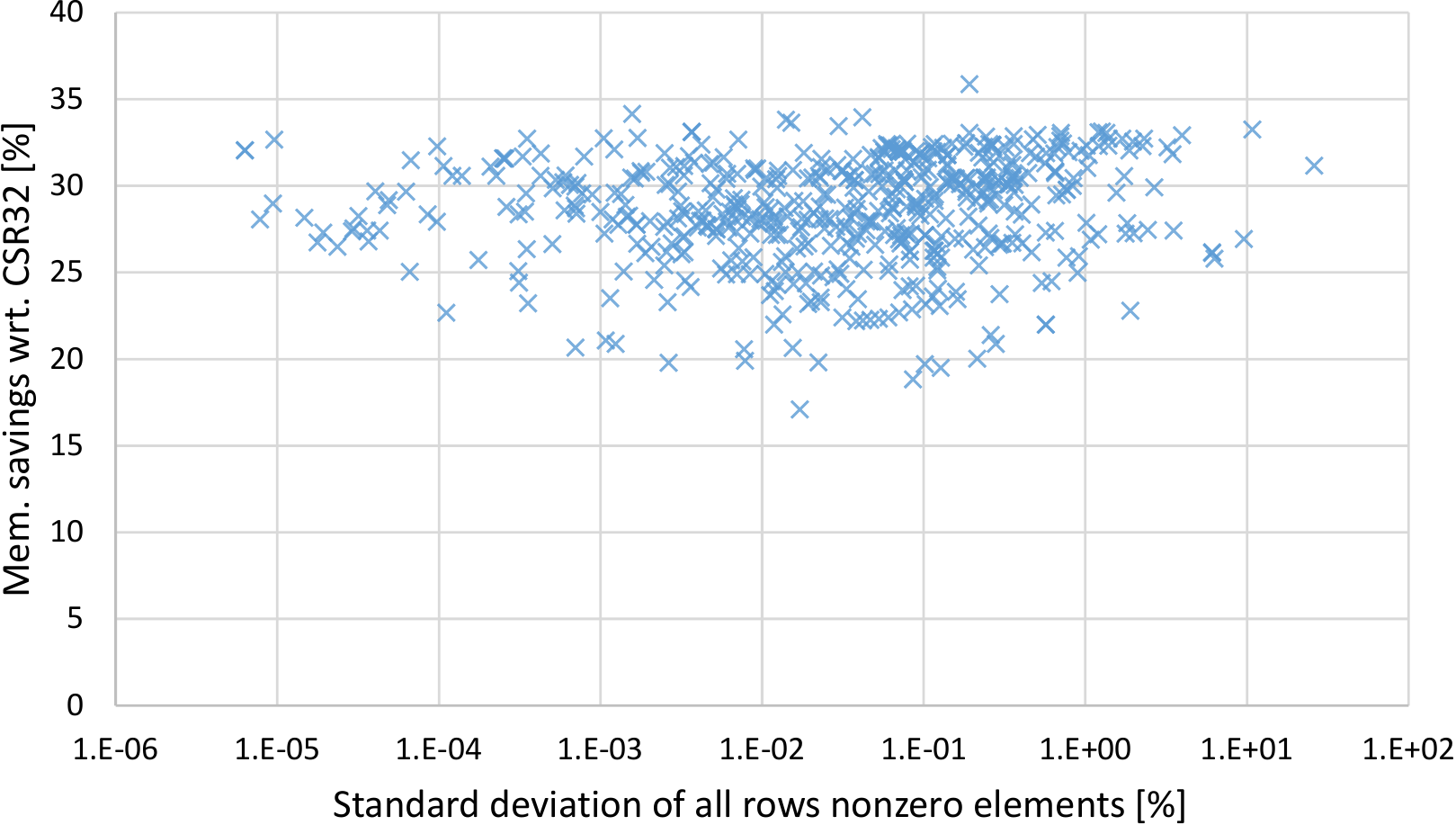}
\end{subfigure}
\hspace{1em}
\begin{subfigure}{.47\textwidth}
	\centering
	\includegraphics[width=\textwidth]{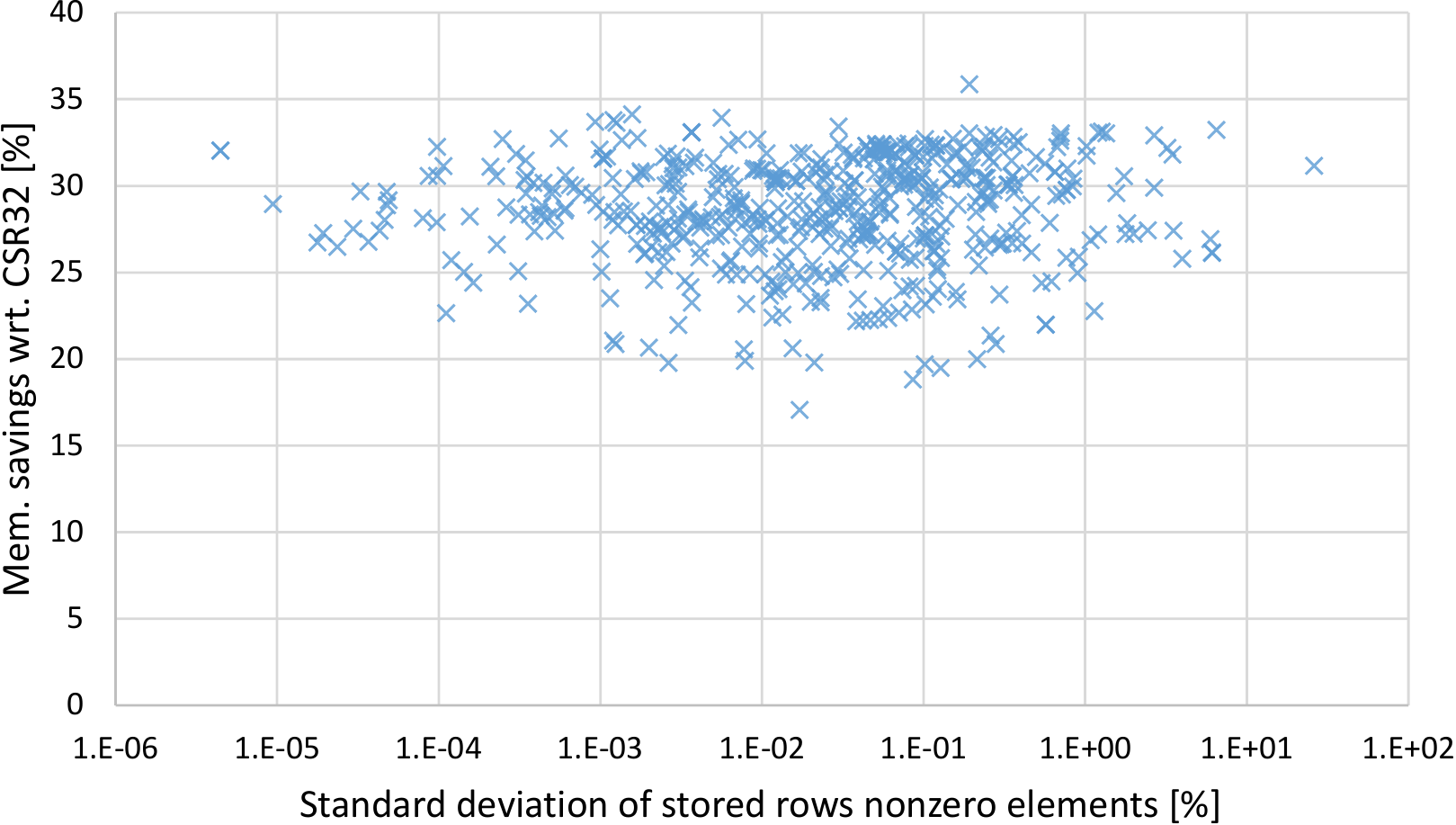}
\end{subfigure}
\caption{Relative memory savings against CSR32 in percents as a function of $\rho$ and $\prnnz$ measured for the tested matrices and double precision considering both all/stored nonzero elements.}
\label{fig:savings:prnnz}
\end{figure}

The measurements for the second and third criteria are presented by the top and bottom parts of \autoref{fig:savings:prnnz}, respectively. Again, we show results only for double precision for the same reason as above. Seemingly (and maybe interestingly), there is no obvious correlation between the memory savings of partitioned matrices against CSR32 and the density of nonzero elements of matrices / uniformity of their distribution across matrix rows.

In the summary, the obtained results support the potential profitability of partitioning of sparse matrices in general.

\subsection{Memory Footprints Compared with Lower Bounds}

\autoref{sec:savings} showed how much memory space we would save if we stored sparse matrices in optimal blocking configurations instead of in CSR32. The last object of our concern within this study was of how much are the memory footprints of the tested matrices higher than their potential minima, i.e., their lower bounds.

We further do not consider compression of the values of matrix nonzero elements, since it is generally worth applying only for special kinds of matrices where nonzero elements contain few unique numbers.
To store $\nnz$ nonzero elements of a matrix $A$ in memory with respect to a precision given by $b$, we thus need $\nnz \times b$ bits to store their values and some additional space to store the information about their structure. The lower bound for the latter for any particular structure of nonzero elements is $1$ bit, since it is sufficient for distinguishing whether or not a matrix has that particular structure. For instance, we can use this bit to indicate whether a matrix is tridiagonal. If it is, the bit would be set and we can store the values of nonzero elements in a dense array; their row and column indices can then be derived from the positions of values in this array. Such an approach can be generally applied for any particular structure of matrix nonzero elements.

In practice, we would likely store in memory also some additional information about a matrix, such as its dimensions or its number of nonzero elements. However, for large matrices such as those from our tested suite, this additional data require a negligible amount of memory, therefore we define a lower bound for a matrix memory footprint simply as $\mmflb(A, b) = \nnz \times b$.

Let
\begin{align}
\Gamma_\boxplus^b(k) &= 
\left(
\frac{ \min \bigl\{ \mmf(A_k, s, h\times w, b) : s\in \schemes_6, h\times w \in \bsizes_{64} \bigr\} }
{\mmflb(A_k,b)}
- 1 \right) \times 100, \\
\Gamma_{\mathrm{CSR32}}^b(k) &= 
\left( \frac{ \mmfcsr(A_k,b) } {\mmflb(A_k,b) } - 1 \right) \times 100. 
\end{align}
$\Gamma_\boxplus^b(k)$ thus expresses of how much percents is the memory footprint of $A_k$ stored in an optimal blocking way higher than its lower bound. For comparison purposes, we define also a corresponding metric for the CSR32 format denoted by $\Gamma_{\mathrm{CSR32}}^b(k)$.

The measured statistics of $\Gamma_\boxplus^b(k)$ and $\Gamma_{\mathrm{CSR32}}^b(k)$ for the tested matrices are shown in \autoref{tab:lb}. Memory footprints of partitioned sparse matrices were obviously much closer to the lower bounds than memory footprints of matrices stored in CSR32; namely, 5 times closer in average and 2 times in worst cases. Moreover, in best cases, partitioned matrices almost reached their lower-bound memory footprints. For instance, in double precision, 7, 26, and 120 matrices out of 563 provided memory footprints up to 1, 2, and 5 percents above their lower bounds, respectively.

\begin{table}[t]
\centering
\small
\begin{tabular}{l*{4}{d{-3}}}
\toprule
& \multicolumn{2}{c}{Single precision} & \multicolumn{2}{c}{Double precision} \\
	\cmidrule(l){2-3} \cmidrule(l){4-5}
Statistics & \multicolumn{1}{c}{Blk.-opt.} & \multicolumn{1}{c}{CSR32} & \multicolumn{1}{c}{Blk.-opt.} & \multicolumn{1}{c}{CSR32} \\ 
	\cmidrule(r){1-1} \cmidrule(l){2-3} \cmidrule(l){4-5}
Minimum &  0.63 & 100.02 &  0.31 & 50.01 \\
Average & 21.85 & 111.03 & 10.93 & 55.51 \\
Maximum & 71.31 & 152.39 & 35.66 & 76.19 \\
\bottomrule
\end{tabular}
\caption{Statistics of $\Gamma_\boxplus^b(k)$ and $\Gamma_{\mathrm{CSR32}}^b(k)$ (in percents) for the tested matrices.}
\label{tab:lb}
\end{table}

\section{Conclusions}

Within this study, we analyzed memory footprints of 563 representative sparse matrices with respect to their partitioning into uniformly sized blocks. We considered different block sizes and different ways of storing blocks in a computer memory. 
The obtained results led us to the following conclusions:
\begin{enumerate}
\item Partitioning of sparse matrices substantially reduce memory footprints of sparse matrices when compared to the most-commonly used storage format CSR32. The average observed memory savings in case of single and double precision were $42.29$ and $28.67$ percents of memory space, respectively. The corresponding worst-case savings were $25.46$ and $17.08$ percents.

\item Partitioning of sparse matrices provides memory footprints much closer to their lower bounds than CSR32. In average, the measured memory footprints for optimal blocking configurations were of only $21.85$ and $10.93$ percents higher than the lower bounds, while the corresponding memory footprints for CSR32 were higher of $111.03$ and $55.51$ percents. Moreover, the memory footprints of matrices most suitable for block processing approach the lower bounds; the amount of memory required for storing information about the structure of nonzero elements of such matrices is relatively negligible.
 
\item For minimization of memory footprints of partitioned sparse matrices in general, we cannot consider only a single format for storing blocks. Instead, we need to choose a format according to the structure of matrix nonzero elements either for all its blocks collectively (min-fixed scheme) all for each block separately (adaptive scheme); the latter approach provides typically lower memory footprints. 

\item For minimization of memory footprints of partitioned sparse matrices in general, we cannot consider only a single block size. However, we can substantially reduce the set of block sizes in the optimization space and still obtain memory footprints close to their optima. In average, the measured memory footprints for the proposed reduced sets of block sizes $\bsizes_{20}$, $\bsizes_{14}$, and $\bsizes_8$ and the min-fixed/adaptive schemes were at most of only $1.51$ percents higher than the optimal values. Even considering square blocks only is thus generally sufficient for minimization of memory footprints of sparse matrices. However, there exist matrices for which the corresponding metrics are significantly higher and are inversely proportional to the number of tested block sizes. One should thus be aware of whether or not his/her matrices fall into this category and if yes, he/she might consider using larger sets of block sizes.

\item The obtained results seem to be consistent across a wide range of real-world matrices arising from multiple applications problems. 

\item There is seemingly no advantage for storing blocks in CSR; without considering this format for blocks, the memory footprints of matrices grow only slightly or not at all. The COO and bitmap formats themselves minimize memory footprints of partitioned sparse matrices, while the dense format is likely the most efficient for related computations.

\item We measured memory savings of partitioned sparse matrices against CSR32 as a function of the following criteria, which are frequently used in the literature: the application problem type, the density of matrix nonzero elements, and the standard deviation of the number of nonzero elements across matrix rows. To our best, we did not find any correlation between the memory savings and these criteria; the blocking approach thus seems reduce memory footprints of sparse matrices in general.
\end{enumerate}

Our findings are encouraging since they show that memory footprints of partitioned sparse matrices can be substantially reduced even when a relatively small block preprocessing optimization space is considered. Whether of not such a reduction pays off in practice depends first of all on the objective one wants to achieve. A big challenge is to improve the performance of memory-bounded sparse matrix operations due to the reduction of matrix memory footprints.
Within our future work, we plan to face this problem at least partially; namely, we will focus on the development of scalable efficient block preprocessing and SpMV algorithms for the min-fixed and adaptive blocking storage schemes, and we will evaluate them experimentally on mainstream HPC architectures.

\section*{Acknowledgements}

The author would like to thank I. \v{S}ime\v{c}ek from the Czech Technical University in Prague for his insightful comments and suggestions that helped to improve the quality of the article.

\section*{References}

\bibliographystyle{elsarticle-harv}

\bibliography{langr}

\end{document}